\documentstyle{amsppt}
\magnification=\magstep1

\hcorrection{.25in}

\TagsOnRight

\NoBlackBoxes

\topmatter
\title
A mixed Hodge structure on a CR manifold
\\
\endtitle
\author
Takao Akahori
\endauthor
\address
Department of Mathematics, Himeji Institute of Technology, 2167
Shosha, Himeji, Hyogo 671-22, Japan
\endaddress

\abstract
The purpose of this work is to propose a mixed Hodge structure over a CR manifold.
As you know, for a CR manifold, Kohn--Rossi cohomology is naturally introduced. 
However, the relation between Kohn--Rossi cohomology and De Rham cohomology is not
so well understood, even in Tanaka's work. We discuss this point.
\endabstract
\endtopmatter

The purpose of this paper is to propose a mixed Hodge structure for CR
manifolds. In our former papers([A1],[A2]), we studied
deformation theory of CR manifolds, and constructed the versal family of
deformations of CR structures. Later, we discussed the smoothness of the
parameter space of this versal family, namely the analogy of Bogomolov
theorem for compact Kaehler manifolds(cf. [A-M1],[A-M2]). However, our
work is still not so definitive. And so, in order to investigate the
versal family more precisely, we have to formulate our notion of a mixed Hodge
structure for general CR manifolds. Our terminology is not the standard one. For
example, algebraic geometrists already introduced the notin of mixed Hodge
structures for algebraic spaces, and even for link structures(see [Dur]). However,
for general CR-structures, the mixed Hodge structure has not been given. And the relation 
with Kohn--Rossi cohomology has never been discussed. In this paper, for general CR
structures, we establish the notion of mixed Hodge structure, which is applicable to deformation
theory of CR-structures. For this purpose, we use the same
method which succeeded in deformation theory. Namely, we start  with
$\Gamma(M,(C\xi)^{\ast}\wedge
\wedge^p(^0T^{\prime})^{\ast}\wedge^q(^0T^{\prime \prime})^{\ast})$(for the notation, see
Sect.1). And consider $F^{p,q}$, which is defined as a subspace in the above
space. And we show that this space admits $d',d^{\prime \prime}$ and
show that our complex recovers the standard Kohn--Rossi cohomology. As you
know, Kohn--Rossi cohomology doesn't depend on the choice of the supplement
vector field. However, our mixed Hodge structure heavily relies on it.
\par
 We should note the recent Rumin's work on the contact geometry (cf.[Ru]).
It is John M Lee who pointed out that our complex quite resembles Rumin's complex 
and also Rumin's work reminds him of our former work on $E_j$-strucutures 
about deformation theory of CR-structures(see [A1],[A2]). However, our point
of view comes from deformation theory (see Sect.2 in this paper), naturally, 
it should lead to our Hodge theory(see [A-M2]).
And as is mentioned already, our complex arises from deformation theory.
So, we think that the view point is different from Rumin's one, and in order
to stress this point, our work is worth publication.

\heading
1. CR-structures and Kohn--Rossi cohomology group 
\endheading

 Let N be a complex manifold. Let $\Omega$ be a strongly pseudo convex domain
with smooth boundary $M=b\Omega$. Then, as is well knwon, over this boundary
$M$, a CR structure is induced from $N$. Namely, we set
$$
^0T^{\prime \prime}=C\otimes TM\cap T^{\prime \prime}N\mid_M.
$$
Then, this $^0T^{\prime \prime}$, a subbundle of $C\otimes TM$ satisfies;
$$
\align
^0T^{\prime \prime}\cap \overline {^0T^{\prime \prime}}=0,dim_C
\frac {C\otimes TM}{^0T^{\prime \prime}+\overline {^0T^{\prime \prime}}}=1 
\tag 1\\
[\Gamma(M,{}^0T^{\prime \prime}),\Gamma(M,{}^0T^{\prime \prime})]\subset
\Gamma(M,{}^0T^{\prime \prime}). \tag 2
\endalign
$$
This pair $(M,{}^0T^{\prime \prime})$ is called a CR-structure. For this
CR structure, we have
$\overline \partial _b$-operator, 
$$
\overline \partial _b ; \Gamma(M,C) \to \Gamma(M,(^0T^{\prime \prime})^{\ast})
$$
by: for $f\in\Gamma(M,C)$, 
$\overline \partial _bf(X)=Xf, X\in{}^0T^{\prime\prime}.$
Like the standard exterior derivative $d$, we have 
$\overline \partial_b$-complex,
$$
\aligned
0\to \Gamma(M,C) @>\overline \partial _b >>\Gamma(M,(^0T^{\prime
\prime})^{\ast})& @>\overline \partial_b >>\Gamma(M,\wedge
^2(^0T^{\prime \prime})^{\ast}) \to \\
\to \Gamma(M,\wedge ^p(^0T^{\prime \prime})^{\ast})&
@>\overline \partial_b >>\Gamma(M,\wedge ^{p+1}(^0T^{\prime
\prime})^{\ast})\to 
\endaligned
$$
 Now we recall the Kohn--Rossi cohomology(cf.[4]). We set a $C^{\infty}$ vector
bundle decomposition
$$
C\otimes TM={}^0T^{\prime}+{}^0T^{\prime \prime}+C\xi,
$$
where $\xi _p \notin{}^0T^{\prime}_p+{}^0T^{\prime \prime}_p$ for every point
$p$ of $M$, and $^0T^{\prime}$ means $\overline {^0T^{\prime \prime}}$. We
fix this and use the notation $T^{\prime}$ for $^0T^{\prime}+C\xi$, that is
to say,
$$
T^{\prime}={}^0T^{\prime}+C\xi, and C\otimes TM={}^0T^{\prime\prime}+T^{\prime}
\tag 1.2.1
$$
We set
$$
C^{p,q}=\Gamma(M,\wedge ^p(T^{\prime})^{\ast}\wedge \wedge ^q(^0T^{\prime
\prime})^{\ast}).
$$
Let $d$ be the exterior differential operator. Then, not like the case
complex manifolds,
$$
d : C^{p,q} \to  C^{p+1,q}+C^{p,q+1}+C^{p+2,q-1}.
$$
This fact was observed by Tanaka(cf.[6]). In fact, for $u$ in $C^{p,q}$,
for $X_i\in T^{\prime}$,\newline $Y_j\in {}^0T^{\prime \prime}$,
$$
\aligned
(du)(X_1,..,X_{p+2},Y_1,..,Y_{q-1})
=&\sum _j(-1)^{j+1}X_ju(X_1,.,\check X_j,.,X_{p+2},Y_1,..,Y_{q-1})\\
+&\sum _i(-1)^{p+2+i+1}Y_iu(X_1,..,X_{p+2},Y_1,.,\check Y_i,.,Y_{q-1})\\
+&\sum _{r<s}(-1)^{r+s}u([X_r,X_s],.,\check X_r,.,\check
X_s,.,Y_1,..,Y_{q-1})\\
+&\sum _{r,s}(-1)^{r+s+p+2}u([X_r,Y_s],.,\check X_r,.,\check Y_s,.,)\\
+&\sum _{r<s}(-1)^{r+s}u([Y_r,Y_s],X_1..,X_{p+2},\check Y_r,.,\check Y_s,..).
\endaligned
$$
Because of $u \in C^{p,q}$, the first line must vanish. Similarly, the second
line, the fourth line, the fifth line vanish. However, the third line may not
vanish. Because the $^0T^{\prime \prime}$ part of $[X_r,X_s]$ might appear
in general CR manifolds.($^0T^{\prime \prime}$ is integrable, but $T^{\prime}$
is not). With this in mind, we set; for $u$ in $C^{p,q}$,
$$
\overline \partial _bu=(du)_{\wedge ^p(T^{\prime})^{\ast}\wedge \wedge
^{q+1}(^0T^{\prime \prime})^{\ast}}.
$$
Here $(du)_{\wedge ^p(T^{\prime})^{\ast}\wedge \wedge
^{q+1}(^0T^{\prime \prime})^{\ast}}$ means the $\wedge
^p(T^{\prime})^{\ast}\wedge \wedge ^{q+1}(^0T^{\prime \prime})^{\ast}$ part 
of $du$. So, we have
$$
C^{p,q} @>\overline \partial _b >>C^{p,q+1}, 
$$
$\overline \partial _b\overline \partial _b=0$. And we have a cohomology
group(Kohn--Rossi cohomology)
$$
H^{p,q}(M)=\frac {Ker  \overline \partial _b\cap C^{p,q}}{\overline \partial
_bC^{p,q-1}}.
$$
This coholomology group has the following mean. Let $U$ be a 
tubular neighborhood of the boundary $M=b\Omega$ in $N$. Then,
$$
H^{p,q}(M)\simeq H^q(U,\Omega ^p_U), 1\le q \le n-1.
$$
Here $H^q(U,\Omega ^p_U)$ means the $\overline \partial $ Dolbeault
cohomology group over $U$. 
 While, over a tubular neighborhood $U$, there is $\partial$ operator ( $U$
is a complex manifold), and we can discuss the mixed Hodge theory over $U$.
The purpose of my work is to study a mixed Hodge theory over CR - manifolds.
Of course, it is possible to introduce a kind of $\partial _b$ operator.
Namely, for $v\in C^{p,q}$, we set
$$
\partial _bv=(dv)_{\wedge ^{p+1}(T^{\prime})^{\ast}\wedge \wedge
^q(^0T^{\prime \prime})^{\ast}}.
$$
Then, we have
$$
C^{p,q} @>\partial _b >>C^{p+1,q}.
$$
However, as is shown, by Tanaka(cf.[6]), unfortunately, in this case, in
general,
$$
\partial _b\partial _b\ne 0
$$
(because $T^{\prime}$ is not integrable). And so $(C^{p,q},d)$ is not even
a double complex.\par

\heading
2. Deformation theory of CR-structures 
\endheading

In this section, we recall deformation theory of CR-structures.\par
$\bold {2.1.  Almost  CR manifolds}$\par
Let $M$ be a $C^{\infty}$ differentiable manifold with real dimension $2n-1$. Let
$E$ be a $C^{\infty}$ subvector bundle of the complexfied tangent bundle $C\otimes TM$
satisfying;
$$
E\cap \overline E=0, dim_C\frac {C\otimes TM}{E+\overline E}=1.
$$
This pair (M,E) is called an almost CR-manifold or an almost CR-structure. Now let
$(M,{}^0T^{\prime \prime})$ be a CR manifold, Then, by using the $C^{\infty}$
vector bundle decomposition (1.2.1), for each point $p$ of $M$, we have a 
homomorphism map from $E$ to $^0T^{\prime \prime}_p$, a composition map of the 
inclusion of $E$ to $C\otimes TM$ and the projection of $C\otimes TM$ to 
$^0T^{\prime \prime}$.\par

$\bold {Definition 2.1.1.}$ Let $(M,{}^0T^{\prime \prime})$ be a CR manifold. $(M,E)$
is called an almost CR manifold with finite distance from $(M,{}^0T^{\prime \prime})$
if and only if the above homomorhism map is an isomorphism map.\par

Then, we have\par

$\bold {Proposition 2.1.2.}$ If $(M,E)$ an almost CR manifold with finite distance
from $(M,^0T^{\prime \prime})$, then there is a $\phi \in 
\Gamma (M,T'\otimes (^0T^{\prime \prime})^{\ast})$ satisfying;
$$
\align
E=&^{\phi}T^{\prime \prime}\\ 
 =&\{X';X'=X+\phi(X),X\in {}^0T^{\prime \prime}\}
\endalign
$$
For the proof, see [A1].\par
 This almost CR structure $(M,{}^{\phi}T^{\prime\prime})$ is a CR structure if and
only if $\phi$ satisfies the following non linear differential equation.
$$
\overline \partial _{T'}^{(1)}\phi + R_2(\phi) + R_3(\phi)=0.
$$
For the notations, see [A1],[A2]. We recall $\overline \partial _{T'}$ - 
cohomology(so called Deformation complex). Namely, we set $\Gamma (M,T')$,
consisting of $T'$ valued global
$C^{\infty}$ sections on $M$, and consider a first order
differential operator
$$
\overline \partial _{T'} ; \Gamma(M,T') \to  \Gamma(M,T'\otimes 
(^0T^{\prime\prime})^{\ast})
$$
defined as follows.\par
 For $Z$ in $\Gamma (M,{}^0T^{\prime\prime})$, and for a $u$ in
$\Gamma (M,T')$,
$$
[fZ,u]=f[Z,u]-u(f)Z.
$$
So we take $T'$-component according to (1.2.1), then we have
$$
[fZ,u]_{T'}=f[Z,u]_{T'}.
$$
With this in mind, we set a first order differential operator by;
$$
\overline \partial _{T'}u(X)=[X,u]_{T'} for  X \in {}^0T^{\prime\prime}.
$$
And like as for scalar valued differential forms, we can introduce 
differential operators
$$
\overline \partial _{T'}^{(i)}:
\Gamma (M,T'\otimes \wedge ^i(^0T^{\prime\prime})^{\ast})
\to \Gamma (M,T'\otimes \wedge ^{i+1}(^0T^{\prime\prime})^{\ast})
$$
and we have a differential complex
$$
\align
0\to \Gamma(M,T')\to \Gamma&(M,T'\otimes (^0T^{\prime \prime})^{\ast})
\to \Gamma(M,T'\otimes \wedge ^2(^0T^{\prime \prime})^{\ast})\to \\               
 \to \Gamma&(M,T'\otimes \wedge ^p(^0T^{\prime \prime})^{\ast})
\to \Gamma(M,T'\otimes \wedge ^{p+1}(^0T^{\prime \prime})^{\ast})\to
\endalign
$$
$\bold {2.2. E_j structure}$\par
Let $(M,{}^0T^{\prime \prime})$ be an orientable CR structure. Then we can introduce 
the Levi-form over $M$. Namely, for each point $p$ of $M$, we set
$$
L(X_p,Y_p)=-\sqrt {-1}[X',\overline {Y'}]_{C\otimes F},
$$
where $X'$(resp. $Y'$) is a $C^{\infty}$ extension of $X_p$(resp. $Y_p$), and
$[X',\overline {Y'}]_{C\otimes F}$ means the projection of 
$[X',\overline {Y'}]$ to $C\otimes F$ 
according to (1.2.1). If this form is positive or negative definite, then our
CR structure $(M,^0T^{\prime \prime})$ is called strongly pseudo convex. Henceforth,
we assume that our CR structure $(M,^0T^{\prime \prime})$ is strongly pseudo convex. 
In this section we introduce a subbundle $E_j$ of $T'\otimes \wedge ^j(^0T^{\prime \prime})$
and show that the cohomology group associated with 
$(\overline \partial^{(i)}_{T'},\Gamma (M,T'\otimes \wedge ^i(^0T^{\prime \prime})^{\ast}))$
can be reduced to the cohomology group associated with 
$(\overline \partial_i,\Gamma (M,E_i))$.\par
 First we set a subspace $\Gamma _i$ of 
$\Gamma (M,T'\otimes \wedge ^i({}^0T^{\prime \prime})^{\ast})$
by;
$$
\Gamma _i=\{ u ;u\in\Gamma (M,{}^0T'\otimes \wedge ^i(^0T^{\prime \prime})^{\ast}),
               (\overline \partial ^{(i)}_{T'}u)_{C\otimes F}=0\},
$$
where $(\overline \partial^{(i)}_{T'}u)_{C\otimes F}$ denotes the projection 
from $T'\otimes \wedge ^i(^0T^{\prime \prime})$ to $C\otimes F\wedge {}^i(^0T^{\prime \prime})$ 
according to (1.2.1). Then we have\par
$\bold {Proposition 2.2.1.}$ There is a subbundle $E_i$ of 
$T'\otimes \wedge ^i(^0T^{\prime \prime})^{\ast}$ to 
$C\otimes F\wedge ^i(^0T^{\prime \prime})^{\ast}$ satisfying;
$$
\Gamma _i = \Gamma (M,E_i).
$$
$\bold {Proof.}$ We show that the mapping 
$(\overline \partial ^{(i)}_{T'}u)_{C\otimes F}$ from 
$\Gamma (M,^0T'\otimes \wedge ^i (^0T^{\prime \prime})^{\ast})$ to 
$\Gamma (M,C\otimes F \wedge ^{i+1}(^0T^{\prime \prime})^{\ast})$ is linear
with respect to $C^{\infty}$ funtions. In fact, for any $C^{\infty}$ function $f$ and
for any section $u$ in 
$\Gamma (M,^0T'\wedge \otimes (^0T^{\prime \prime})^{\ast})$, we obtain
$$
\aligned
(\overline \partial ^{(i)}_{T'}fu&)_{C\otimes F}(X_1,.,X_{i+1})\\ 
 =&(\sum_j(-1)^j[X_j,fu(X_1,.,\check X_j,.,X_{i+1})]_{T'}\\
 +&\sum_{\alpha < \beta}(-1)^{\alpha + \beta}fu([X_{\alpha},
X_{\beta}],X_1,.,\check X_{\alpha},.,\check X_{\beta},.,X_{i+1}))_{C\otimes F}\\
 =&\sum _j(-1)^j[X_j,fu(X_1,.,\check X_j,.,X_{i+1})]_{C\otimes F}.
\endaligned
$$
As $u$ is an element in $\Gamma (M,^0T'\otimes \wedge ^i(^0T^{\prime \prime})^{\ast})$,
this leads
$$
\aligned
(\overline \partial ^{(i)}_{T'}fu&)_{C\otimes F}(X_1,.,X_{i+1})\\
 =&\sum _j(-1)^j[X_j,fu(X_1,.,\check X_j,.,X{i+1})]_{C\otimes F}\\
 =&f\sum _j(-1)^j[X_j,u(X_1,.,\check X_j,.,X_{i+1})]_{C\otimes F}\\
 =&f(\overline \partial ^{(i)}_{T'}u)_{C\otimes F}(X_1,.,X_{i+1}) for all \Gamma (M,^0T^{\prime \prime}).
\endaligned
$$
Therefore for each point $p$ of $M$, we can deffine a vector space $E_{i,p}$ by;
$$
E_{i,p}=\{u_p;u\in ^0T'\otimes \wedge ^i({}^0T^{\prime
\prime})^{\ast},(\overline \partial ^{(i)}_{T'}u)_{C\otimes F}(p)=0\}.
$$
First, we prove that $dim_CE_{i,p}$ is independent of $p$. To do this, it is enough
to show that the following vector bundle sequence is exact.
$$
\align
^0T'\otimes \wedge ^i(^0T^{\prime \prime})^{\ast}&\to 
C\otimes F\otimes \wedge^{i+1}(^0T^{\prime \prime})^{\ast}\to 0 \tag 2.3.1\\
                               u&\to (\overline \partial ^{(1)}_{T'}u)_{C\otimes F}
\endalign
$$
For each point $p$ in $M$, we take a system of moving frame $\{e_j\}_{j=1,.,n-1}$
of $^0T^{\prime \prime}$ in a neighborhood of $p$ and $F$ such that
$$
[e_i,\overline e_j]_{C\otimes F}=\sqrt {-1}\delta_{i,j}F.
$$
Then for any $F\otimes e^{\ast}_{k_1}\wedge e^{\ast}_{k_2}\wedge...\wedge e^{\ast}_{k_{i+1}}
(k_1<k_2<...<k_{i+1})$, we set
$$
u=\overline e_{k_1}\otimes e^{\ast}_{k_2}\wedge ...\wedge e^{\ast}_{k_{i+1}}.
$$
By a simple calculation, we obtain
$$
(\overline \partial ^{(1)}_{T'}u)_{C\otimes F}=F\otimes e^{\ast}_{k_1}\wedge ...\wedge e^{\ast}_{k_{i+1}}.
$$
Therefore the above mapping $(\overline \partial^{(i)}_{T'}u)_{C\otimes F}$ is
surjective. Thus the sequence (2.2.1) is exact and $dim_CE_{i,p}$ is indepent of $p$.
And so, $\cup _{p\in M}E_{i,p}$ is a vector bundle on $M$. By $E_j$, we write this vector 
bundle. Then, we have
$$
\Gamma _i=\Gamma (M,E_i). \bold{Q.E.D}
$$ 
We see this vector bundle more precisely.\par
$\bold{Theorem 2.2.2.}$ $E_0=0$. And there is a following differential subcomplex.
$$
\align
0 \to  \Gamma (M,E)_1) @>\overline \partial_1 >> 
\Gamma& (M,E_2) @>\overline \partial_2 >>\Gamma (M,E_3)\\
@>\overline \partial_{i-1}>> \Gamma& (M,E_i) @>\overline \partial_i >>
\Gamma (M,E_{i+1})
\endalign
$$
where $\overline \partial_i$ means the restriction of 
$\overline \partial^{(i)}_{T'}$ to $\Gamma (M,E_i)$.\par
$\bold{Proof}.$ By the definition of $E_0$, we have
$$
E_0=\{u;u\in{}^0T', [u,X]_{C\otimes F}=0 for all X\in{} ^0T^{\prime
\prime}\}.
$$
Since the Levi form is non-degenerate, we obtain
$$
E_0=0.
$$
Next we prove
$$
\overline \partial^{(i)}_{T'}\Gamma(M,E_i)\subset \Gamma (M,E_{i+1}).
$$
In fact for all $u\in \Gamma (M,E_i)$, it follows
$$
\overline \partial^{(i)}_{T'}u\in \Gamma (M,^0T'\otimes \wedge ^{i+1}(^0T^{\prime \prime})^{\ast})
$$
from the relation $(\overline \partial^{(i)}_{T'}u)_{C\otimes F}=0$(definition of 
$E_i$). And
$$
(\overline \partial^{(i+1)}_{T'}(\overline \partial^{(i)}_{T'}u))_{C\otimes F}=0 
(from \overline \partial^{(i+1)}_{T'}\overline \partial ^{(i)}_{T'}=0).
$$
These considerations leads to that;for any $u\in \Gamma (M,E_i)$, we have
$$
\overline \partial^{(i)}_{T'}u\in \Gamma(M,E_{i+1}).
$$
This completes the proof. $\bold{Q.E.D.}$\par
 As for $E_1$, we prove that $\Gamma (M,E_1)$ has sufficiently many sections. We
impose the Levi-metric over M and we form the adjoint operator 
$\overline \partial ^{\ast}_{T'}$ of $\overline \partial _{T'}$. We set the
Laplacian
$$
\square _{T'}=\overline \partial _{T'}\overline \partial ^{\ast}_{T'}
              + \overline \partial ^{\ast}_{T'}\overline \partial _{T'}.
$$
As is usual in the theory of harmonic forms, we obtain the harmonic space
$\bold {H^{(1)}_{T'}}$ in $\Gamma (M,T'\otimes (^0T^{\prime \prime})^{\ast})$
(we assume that $dim_RM\ge 7$). We introduce a differential operator $\Cal L$
from $\Gamma (M,T'\otimes (^0T^{\prime \prime})^{\ast})$ to
$\Gamma (M,^0T'\otimes (^0T^{\prime \prime})^{\ast})$ as folows; for each $\phi$,
$\Gamma (M,T'\otimes (^0T^{\prime \prime})^{\ast})$, we put
$$
\Cal L\phi = \phi - \overline \partial _{T'}\theta _{\phi},
$$
where $\theta _{\phi}$ is an element of $\Gamma (M,{}^0T')$ defined by;
$$
[X,\theta _{\phi}]_{C\otimes F}=\phi (X)_{C\otimes F} for any X \in \Gamma
(M,{}^0T^{\prime \prime}).
$$
Then we have the following theorem.\par
$\bold {Theorem 2.2.3.}$ The mapping $\Cal L\mid _{\bold {H^{(1)}_{T'}}}$, being 
restricted to $\bold {H^{(1)}_{T'}}$, is injective and 
$$
\Cal H  \subset  \Gamma (M,E_1)
$$
holds, where $\Cal H$ denotes $\Cal L(\bold {H^{(1)}_{T'}})$.\par

$\bold{Proof.}$ Injectivity is obvious. So it suffices to show
$$
\Cal H   \subset   \Gamma (M,E_1).
$$
That is to say, for all $\phi \in \bold {H^{(1)}_{T'}}$,
$$
\Cal L\phi \in \Gamma (M,{}^0T'\otimes (^0T^{\prime \prime})^{\ast})
$$
and
$$
(\overline {\partial} ^{(1)}_{T'}(\Cal L\phi ))_{C\otimes F})=0.
$$
By the definition of $\Cal L$, we have that; for $\phi$ in 
$\Gamma (M,T'\otimes (^0T^{\prime \prime})^{\ast})$
$$
\Cal L\phi  in  \Gamma (M,^0T'\otimes (^0T^{\prime \prime})^{\ast})
$$
and for $\phi$ in $Ker \overline \partial ^{(1)}_{T'}$, obviously,
$$
\Cal L\phi  in  Ker \overline \partial ^{(1)}_{T'}.
$$
Therefore our theorem follows.  $\bold{Q.E.D.}$\par
 
 Especially by Theorem 2.2.3, we have that the injection  $i ;\Cal H
\hookrightarrow Ker \overline \partial ^{(1)}_{T'}$ induces the surjective
map
$$
\Cal H  \to  Ker \overline \partial ^{(1)}_{T'}/Im \overline \partial _{T'}.
$$
 As for $E_i (2\le i \le n-1)$, we have the following theorem.\par
$\bold {Theorem 2.2.4.}$ The injection induces the isomorphism map
$$
\tau: Ker \overline \partial _i/Im \overline \partial _{i-1}\to 
         Ker \overline \partial ^{(i)}_{T'}/Im \overline \partial ^{(i-1)}_{T'},
$$
where $2\le i \le n-1$.\par
$\bold {Proof}.$ First, we prove that the above map is surjective. For this purpose,
it suffices to prove that for all $\phi$ in $Ker \overline \partial ^{(i)}_{T'}$,
there is an element $\theta _{\phi}$ of 
$\Gamma (M,^0T'\otimes \wedge ^{i-1}(^0T^{\prime \prime})^{\ast})$ satisfying;
$$
\phi - \overline \partial ^{(i-1)}_{T'}\theta _{\phi} \in \Gamma (M,E_i).
$$
$\phi$ being in $Ker \overline \partial ^{(i)}_{T'}$, we obtain
$$
\phi - \overline \partial ^{(i-1)}_{T'}\theta _{\phi} \in \Gamma (M,E_i).
$$
Therefore the surjectivity is proved. Next, we prove that the above map is
in injective. If $\psi$ is in $Ker \overline \partial _i$ satisfying 
$\psi = \overline \partial ^{(i-1)}_{T'}\phi$, there is an element $\theta _{\phi}$
in $\Gamma (M,^0T'\otimes \wedge ^{(i-2)}(^0T^{\prime \prime})^{\ast})$ such
that
$$
(\psi)_{C\otimes F}=(\overline \partial ^{(i-2)}_{T'}\theta _{\phi})_{C\otimes F}
$$
by the same argument as above(we assume $i\ge 2$). Of course 
$\phi - \overline \partial ^{(i-2)}_{T'}\theta _{\phi}$ is in $\Gamma (M,E_i)$.
And
$$
\psi=\overline \partial ^{(i-1)}_{T'}(\phi-\overline \partial ^{(i-2)}_{T'}\theta _{\phi})
$$
holds. So we have our theorem. $\bold {Q.E.D.}$
 
\heading
3. $F^{p,q}$ complex and Mixed Hodge structure 
\endheading

 We start with
$$
D^{p,q}(M)=\Gamma(M,(C\xi)^{\ast}\wedge \wedge
^{p-1}(^0T^{\prime})^{\ast}\wedge \wedge ^q(^0T^{\prime \prime})^{\ast}),if 
p\ge 1,
$$
$$
D^{0,q}(M)=0.
$$
And consider
$$
F^{p,q}(M)=\{u: u\in D^{p,q}(M), (du)_{\wedge ^p(^0T^{\prime})^{\ast}\wedge
\wedge ^{q+1}(^0T^{\prime \prime})^{\ast}}=0\},
$$
where $(du)_{\wedge ^p(^0T^{\prime})^{\ast}\wedge
\wedge ^{q+1}(^0T^{\prime \prime})^{\ast}}$ means the $\wedge
^p(^0T^{\prime})^{\ast}\wedge
\wedge ^{q+1}(^0T^{\prime \prime})^{\ast}$ part of $du$.
We see this $F^{p,q}$ more precisely. For this, we compute
$(du)_{\wedge ^p(^0T^{\prime})^{\ast}\wedge
\wedge ^{q+1}(^0T^{\prime \prime})^{\ast}}$. For $X_i \in {}^0T'$,$Y_j \in
{}^0T^{\prime \prime}$,
$$
\aligned
(&du)(X_1,..,X_p,Y_1,..,Y_{q+1})\\
=&\sum _j(-1)^{i+1}X_ju(X_1,.,\check X_i,.,X_p,Y_1,..,Y_{q+1})\\
+&\sum _j(-1)^{p+j+1}Y_iu(X_1,..,X_p,Y_1,.,\check Y_j,.,Y_{q+1})\\ 
+&\sum _{r<s}(-1)^{r+s}u([X_r,Y_s],X_1,.,\check X_r,.,\check
X_s,.,Y_1,..,Y_{q+1})\\
+&\sum _{r,s}(-1)^{r+p+s}u([X_r,Y_s],X_1,.,\check
X_r,.,X_p,Y_1,.,\check Y_s,.,Y_{q+1})\\ +&\sum
_{r<s}u([Y_r,Y_s],X_1,..,X_p,Y_1,.,\check Y_r,.,\check Y_s,.,Y_{q+1})
.\endaligned
$$
So the condition $(du)_{\wedge ^p(^0T^{\prime})^{\ast}\wedge
\wedge ^{q+1}(^0T^{\prime \prime})^{\ast}}=0$ becomes
$$
d\theta \wedge u=0,
$$
where $\theta$ is a real one form defined by;
$$
\theta(\xi)=1, \theta \mid _{^0T^{\prime}+{}^0T^{\prime\prime}}=0.
$$
We set
$$
F^k=\sum _{p+q=k}F^{p,q}.
$$
 Now we introduce
$d^{\prime},d^{\prime \prime}$ operators by;\par
for $u$ in $F^{p,q}$, 
$$
d^{\prime}u=(du)_{(C\xi)^{\ast}\wedge \wedge ^p(^0T^{\prime})^{\ast}\wedge
\wedge ^q(^0T^{\prime \prime})^{\ast}},
$$
\par
for $u$ in $F^{p,q}$,
$$
d^{\prime \prime}u=(du)_{(C\xi)^{\ast}\wedge \wedge
^{p-1}(^0T^{\prime})^{\ast}\wedge \wedge ^{q+1}(^0T^{\prime \prime})^{\ast}}.
$$
Then, our theorem is;\par
$\bold {Theorem 3.1}$.
$$
\aligned
d&^{\prime}F^{p,q}\subset F^{p+1,q},\\
d&^{\prime \prime}F^{p,q}\subset F^{p,q+1},\\
d&F^k\subset F^{k+1}\\
d&^{\prime}d^{\prime}=d^{\prime \prime}d^{\prime \prime}=0,\\
d&^{\prime}d^{\prime \prime}+d^{\prime \prime}d^{\prime}=0\\
d&=d^{\prime}+d^{\prime \prime}.
\endaligned
$$
Namely our $(F^k,d^{\prime},d^{\prime \prime})$ is a double complex.
Then, we have three cohomology groups which were observed in [A-M2].
The first one is
$$
\frac {Ker d^{\prime\prime} \cap  F^{p,q}}{d^{\prime\prime}F^{p,q-1}}.
$$
The second one is
$$
\frac {Ker d^{\prime} \cap  Ker d^{\prime\prime} \cap  F^{p,q}}
{Ker d^{\prime} \cap d^{\prime\prime}F^{p,q-1}}
$$
because of
$$
d^{\prime\prime}(Ker d^{\prime}) \subset Ker d^{\prime}.
$$
The third one is
$$
\frac { d^{\prime} \cap  Ker d^{\prime\prime} \cap  F^{p,q}}
{d^{\prime}F^{p-1,q} \cap  d^{\prime\prime}F^{p,q-1}}
$$
because of
$$
d^{\prime\prime}d^{\prime} = d^{\prime}(-d^{\prime\prime}).
$$
In fact, we have the following theorem.
\par
$\bold {Theorem 3.2}$. If $(M,^0T^{\prime \prime})$ is strongly pseudo
convex with $dim_RM=2n-1$, as for $d^{\prime \prime}$ - operator,
\par if $p+q=n$,
$$
Ker d^{\prime \prime} \to \frac {Ker \overline \partial _b\cap
C^{p,q}}{\overline \partial _bC^{p,q-1}} \to 0,
$$
\par
if $p+q\ge n+1$,
$$
\frac {Ker d^{\prime \prime}}{Im d^{\prime \prime}}\simeq
\frac {Ker \overline \partial _b\cap C^{p,q}}
{\overline \partial _bC^{p,q-1}},
$$
and as for $d$ -operator,
\par if $k=n$,
$$
Ker d \cap F^k \to \frac {Ker d}{Im d} \to 0,
$$
\par if $k\ge n+1$,
$$
\frac {Ker d\cap F^k}{dF^{k-1}}\simeq \frac {Ker d}{Im d}.
$$
So our double complex $(F^k,d^{\prime},d^{\prime\prime})$
recovers the standard  Kohn--Rossi cohomology of  degree
$\ge n+1$, and so we can discuss a mixed Hodge theory over CR manifolds.
The proof of Theorem 3.2 is very like in the proof for $E_j$ bundles(cf.[2]). 
We consider a bundle map from $(C\xi)^{\ast}\wedge \wedge
^{l-1}(^0T^{\prime})^{\ast}\wedge \wedge ^{s-1}(^0T^{\prime \prime})^{\ast}$
to $\wedge ^l(^0T^{\prime})^{\ast}\wedge
\wedge ^s(^0T^{\prime \prime})^{\ast}$, defined by :
$$
(C\xi)^{\ast}\wedge \wedge ^{l-1}(^0T^{\prime})^{\ast}\wedge \wedge
^{s-1}(^0T^{\prime \prime})^{\ast} \to \wedge ^l(^0T^{\prime})^{\ast}\wedge
\wedge ^s(^0T^{\prime \prime})^{\ast}
$$
$$
u \mapsto (du)_{\wedge ^l(^0T^{\prime})^{\ast}\wedge
\wedge ^s(^0T^{\prime \prime})^{\ast}}.
$$
The key lemma is as follows.\par
$\bold {Lemma 3.3}$.
The above map 
$$
(C\xi)^{\ast}\wedge \wedge ^{l-1}(^0T^{\prime})^{\ast}\wedge \wedge
^{s-1}(^0T^{\prime \prime})^{\ast} \to \wedge ^l(^0T^{\prime})^{\ast}\wedge
\wedge ^s(^0T^{\prime \prime})^{\ast},
$$
is surjective if $l+s\ge n$.\par
$\bold {Proof}$. For $D^{p,q}(M)=\Gamma(M,(C\xi)^{\ast}\wedge \wedge
^{p-1}(^0T^{\prime})^{\ast}\wedge \wedge ^q(^0T^{\prime\prime})^{\ast})$, we
have an operator
$$
L: D^{p,q}(M)\to D^{p+1,q+1}(M)
$$
by $u \mapsto d\theta \wedge u$. \par
Then, just like the case hermitian manifolds, we obtain the adjoint operator
of $L$, $\Lambda$, and we can introduce the notion of {\it primitive forms}.
Namely, for $u$ in $D^{p,q}(M)$, $u$ is called a primitive form if and only
if
$$
\Lambda u=0.
$$
For primitive forms, we have two fact.\par
$\bold {Fact.1}$. For $v$ in $D^{p,q}(M)$, $v$ is primitive if and only if
$$
\align
L^rv=0, where r=&max(n-1-(p-1+q),0)\\
              =&max(n-p-q,0).
\endalign
$$
\par
$\bold {Fact.2}$. For $u$ in $D^{p,q}(M)$, $u$ has the following unique
decomposition.
$$
u=u_0+Lu_1+...+L^ku_k, where k=[\frac {p-1+q}2].
$$
The proof is the same as in hermitian manifolds. So we omit this.\par
 Now we see Lemma 3.3. In order to show Lemma 3.3, it suffices to show that;
if $l+s\ge n$, then the following map is surjective.
$$
L: (C\xi)^{\ast}\wedge \wedge ^{l-1}(^0T^{\prime})^{\ast}\wedge \wedge
^{s-1}(^0T^{\prime \prime})^{\ast} \to 
(C\xi)^{\ast}\wedge \wedge ^l(^0T^{\prime})^{\ast}\wedge \wedge
^s(^0T^{\prime \prime})^{\ast}
$$
by $u \mapsto d\theta\wedge u$, because of the computation of
$$
(du)_{\wedge ^l(^0T^{\prime})^{\ast}\wedge \wedge ^s(^0T^{\prime
\prime})^{\ast}}.
$$
While for $v$ in $(C\xi)^{\ast}\wedge \wedge ^l(^0T^{\prime})^{\ast}\wedge
\wedge ^s(^0T^{\prime \prime})^{\ast}$ $(l+s\ge n)$, $v$ is primitive if and
only if $v=0$(by Fact 1). And by Fact 2, we have : for $u$ in
$(C\xi)^{\ast}\wedge \wedge ^l(^0T^{\prime})^{\ast}\wedge
\wedge ^s(^0T^{\prime \prime})^{\ast}$,
$$
u=u_0+Lu_1+....+L^ku_k, where k =[\frac {l+s}2],
$$
and $u_i$ are primitive. So in our case, $u_0$ must be zero. Hence
$$
\align
u=&Lu_1+...+L^ku_k\\
 =&L(u_1+...+L^{k-1}u_k).
\endalign
$$
Therefore we have the surjectivity.

\heading
4. Estimates 
\endheading

By the same method as in [A1], we show the following a priori estimate.
For this, we put the Levi metric on $F^{p,q}$ and consider the adjoint 
operator of $d^{\prime \prime}(resp.  d^{\prime}),d^{\prime \prime \ast}(resp. 
d^{\prime \ast})$.\par
$\bold {Theorem 4.1}$. If $2n+1\ge k=p+q\ge n+1$ and $n-1\ge p,q \ge 2$, then our complexes $(F^k,d)$,
$(F^{p,q},d^{\prime})$, $(F^{p,q},d^{\prime \prime})$ are subelliptic.\par
 
 In order to prove our theorem, we must prepare several facts. On
$\Gamma(M,(C\xi)^{\ast}\wedge\wedge ^{p-1}(^0T^{\prime})^{\ast}\wedge \wedge
^q(^0T^{\prime\prime})^{\ast})$, we we can introduce the formal adjoint
operators, $\delta ^{\prime\prime}$ of $d^{\prime\prime}$, and $\delta
^{\prime}$ of $d^{\prime}$ as in for $\Gamma(M,\wedge
^p(^0T^{\prime})^{\ast}\wedge
\wedge ^q(^0T^{\prime\prime})^{\ast})$(see [T]). By using this operator, we
compute the adjoint operator of $d^{\prime\prime}$ for $F^{p,q}$. 
\par
$\bold {Lemma 4.2}$. The adjoint operator of $d^{\prime\prime}$ in $F^{p,q}$
becomes ; for $u$ in $F^{p,q}$,
$$
d^{\prime\prime\ast}u=\delta ^{\prime\prime}u - \frac 1{p+q-n+1}\Lambda L \delta
^{\prime\prime}u.
$$

$\bold {Proof}$. First, we show that; for $u$ in $F^{p,q}$,
$$
\delta ^{\prime\prime}u-\frac 1{p+q-n+1}\Lambda L\delta ^{\prime\prime}u
\in  F^{p,q-1}
$$ 
For the proof, it suffices to show ; for $v$ in $F^{p,q-1}$,
$u$ in $F^{p,q}$,
$$
\aligned
(&d^{\prime\prime}v,u)\\
=(&v,(\delta ^{\prime\prime}-\frac 1{p+q-n+1}\Lambda L\delta ^{\prime\prime})u)
\endaligned
$$
$\bold {Q.E.D.}$.
\par
 Now we establish an a priori estimate for $\|d^{\prime\prime}u\|^2
+\|d^{\prime\prime\ast}u\|^2$ for $u$ in $F^{p,q}$. For this, we recall several
Kaehler identities on $\Gamma (M,(C\zeta)^{\ast}\wedge\wedge 
^p(^0T^{\prime})^{\ast}\wedge\wedge ^q(^0T^{\prime\prime})^{\ast})$.
\par
$\bold {Kaehler identities}$.
$$
\aligned
[&\Lambda,d^{\prime\prime}]=-\sqrt {-1}\delta ^{\prime}\\
[&\Lambda,d^{\prime}]=\sqrt {-1}\delta ^{\prime\prime}\\
[&\delta ^{\prime\prime},L]=\sqrt {-1}d^{\prime}\\
[&\delta ^{\prime},L]=-\sqrt {-1}d^{\prime\prime}\\
\delta &^{\prime}d^{\prime\prime}+d^{\prime\prime}\delta ^{\prime}=0\\
\delta &^{\prime\prime}d^{\prime}+d^{\prime}\delta ^{\prime\prime}=0\\
[&L,\Lambda]=k-(n-1), where k=p+q
\endaligned
$$
With these equalities, we establish an a priori estimte.
For $u \in F^{p,q}$,
$$
\aligned
\|\Lambda d^{\prime}u\|^2=&(\Lambda d^{\prime}u,\Lambda d^{\prime}u)\\
                         =&(L\Lambda d^{\prime}u,d^{\prime}u),\\
\|L\delta ^{\prime\prime}u\|^2=&\|d^{\prime}u\|^2
\endaligned
$$
By the way, because of $(L\Lambda -\Lambda L)v=(k-(n-1))v$ for
$k$-form $v$,
$$
L\Lambda d^{\prime}u=(k-n+1)d^{\prime}u.
$$
So,
$$
\aligned
\|d^{\prime \prime \ast}u\|^2=&\|\delta ^{\prime \prime}u-\frac
1{p+q-n+1}\Lambda L\delta ^{\prime \prime}u\|^2\\
=&\|\delta ^{\prime \prime}u\|^2+\|\frac 1{p+q-n+1}\Lambda L\delta
^{\prime \prime}u\|^2-2Re(\delta ^{\prime \prime}u,\frac
1{p+q-n+1}\Lambda L\delta ^{\prime\prime}u)\\
=&\|\delta ^{\prime\prime}u\|^2+\frac 1{(p+q-n+1)^2}\|\Lambda L\delta
^{\prime\prime}u\|^2-\frac 2{p+q-n+1}(L\delta ^{\prime\prime}u,
L\delta ^{\prime\prime}u)\\
=&\|\delta ^{\prime\prime}u\|^2+\frac 1{p+q-n+1}\|d^{\prime}u\|^2
-\frac 2{p+q-n+1}\|d^{\prime}u\|^2\\
=&\|\delta ^{\prime\prime}u\|^2-\frac 1{p+q-n+1}\|d^{\prime}u\|^2
\endaligned
\tag 4.1
$$
Hence
$$
\align
\|&d^{\prime\prime}u\|^2+\|d^{\prime \prime \ast}u\|^2\\
=&\|d^{\prime\prime}u\|^2+\|\delta ^{\prime\prime}u\|^2-
\frac 1{p+q-n+1}\|d'u\|^2\\
\ge &\|d^{\prime\prime}u\|^2+\|\delta ^{\prime\prime}u\|^2-
\frac 1{p+q-n+1}\{\|d'u\|^2+\|\delta 'u\|^2\}\\
\ge &\sum_{I,J}\{(\sum _{k\notin J}\|e_ku_{I,J}\|^2
+\sum_{j\in J}\|\overline e_ju_{I,J}\|^2)\\
-&\frac 1{p+q-n+1}(\sum_{l\notin I}\|\overline e_lu_{I,J}\|^2
+\sum _{i\in I}\|e_iu_{I,J}\|^2)\} ,where \mid I\mid =p-1, \mid J \mid =q\\
\ge & \frac {p+q-n}{p+q-n+1}\sum _{I,J}\{\sum _{k\notin J} \|e_ku_{I,J}\|^2
+\sum_{j\in J}\|\overline e_ju_{I,J}\|^2\}\\
+& \frac 1{p+q-n+1}\{\sum_{k\notin J} \|e_ku_{I,J}\|^2
+\sum_{j\in J}\|\overline e_ju_{I,J}\|^2\\
-&\sum_{l\notin I}\|\overline e_lu_{I,J}\|^2-\sum_{i\in I}\|e_iu_{I,J}\|^2\}
\endalign
$$
While
$$
\align
J^c&=I\cap J^c+I^c\cap J^c , J=I\cap J+I^c\cap J\\
I^c&=I^c\cap J^c+I^c\cap J , I=I\cap J^c+I\cap J
\endalign
$$
And so

$$
\sum _{k\notin J}\|e_ku_{I,J}\|^2
=\sum _{k\in I\cap J^c}\|e_ku_{I,J}\|^2+\sum_{k\in I^c\cap J^c}
\|e_ku_{I,J}\|^2
$$
$$
\sum _{j\in J}\|\overline e_ju_{I,J}\|^2
=\sum _{j\in I\cap J^c}\|\overline e_ku_{I,J}\|^2+\sum_{j\in I^c\cap J^c}
\|\overline e_ju_{I,J}\|^2
$$
$$
\sum _{l\notin I}\|\overline e_lu_{I,J}\|^2
=\sum _{l\in I^c\cap J^c}\|\overline e_lu_{I,J}\|^2+\sum_{l\in I^c\cap J}
\|\overline e_lu_{I,J}\|^2
$$
$$
\sum _{i\in I}\|e_ku_{I,J}\|^2
=\sum _{i\in I\cap J^c}\|e_iu_{I,J}\|^2+\sum_{i\in I\cap J}
\|e_iu_{I,J}\|^2
$$
By using this, the above becomes
$$
\aligned
\ge & \frac {p+q-n}{p+q-n+1}\sum _{I,J}\{\sum _{k\notin J} \|e_ku_{I,J}\|^2
+\sum_{j\in J}\|\overline e_ju_{I,J}\|^2\}\\
+& \frac 1{p+q-n+1}\{\sum_{l\in I\cap J}\|\overline e_lu_{I,J}\|^2
-\sum_{l\in I\cap J}\|e_lu_{I,J}\|^2\\
+&\sum _{m\in I^c\cap J^c}\|e_mu_{I,J}\|^2
-\sum _{m\in I^c\cap J^c}\|\overline e_mu_{I,J}\|^2\}
\endaligned
\tag 4.2
$$
While
$$
\aligned
\mid I \mid + \mid J \mid - \mid I\cap J\mid=&\mid I\cup J\mid \\
                                            =&n-1 - \mid I^c\cap J^c\mid
\endaligned
$$
So
$$
\aligned
\mid I\cap J\mid=&p-1+q-n+1+\mid I^c\cap J^c\mid \\
=& p+q-n +\mid I^c\cap J^c\mid.
\endaligned
$$
Especially,
$$
\mid I^c\cap J^c\mid \le \mid I\cap J\mid (by p+q-n\ge 1).
$$
And
$$
\mid I\cap J\mid\le (p+q-n)\mid J^c\mid + \mid I^c \cap J^c\mid
$$
These mean that there is an injective map
$$
\kappa;  I\cap J  \mapsto  (p+q-n)J^c\cup (I^c\cap J^c)
$$
satisfying
$$
\kappa (I\cap J)  \supset  I^c\cap J^c,
$$
where $(p+q-n)J^c$ means the disjoint $(p+q-n)$'s union of $J^c$.
Hence
$$
\aligned
\|&\overline e_lu_{I,J}\|^2+\|e_{\kappa (l)}u_{I,J}\|^2\\
\sim & \|e_lu_{I,J}\|^2+\|\overline e_{\kappa (l)}u_{I,J}\|^2
\endaligned
$$
In the case $\mid J^c\mid =1$,
$$
\|d^{\prime\prime}u\|^2+\|d^{\prime\prime \ast}u\|^2 \gtrapprox
\|d'u\|^2.
$$
By (4.1) with (4.2), we have
$$
\|d^{\prime\prime}u\|^2+\|d^{\prime\prime \ast}u\|^2+\|u\|^2
\gtrapprox \|d^{\prime\prime}u\|^2+\|\delta ^{\prime\prime}u\|^2+\|u\|^2
\gtrapprox \|u\|^{\prime 2}
$$
In the case $\mid J^c\mid >1 $,
we directly have 
$$
\|d^{\prime\prime}u\|^2+\|d^{\prime\prime\ast}u\|^2+\|u\|^2 \gtrapprox
\|u\|^{\prime 2}.
$$

\heading
5. Finiteness
\endheading

In Sect.3, we introduced three groups. In this section, we show
\par
$\bold {Theorem 5.1}$. If $2n+1\ge p+q\ge n+1$ and $n-2\ge p,q\ge 2$,
$$
\frac
{Ker d^{\prime} \cap  Ker d^{\prime\prime} \cap  F^{p,q}}
{Ker d^{\prime} \cap  d^{\prime\prime}F^{p,q-1}}
\backsimeq
\bold {H_{d^{\prime\prime}}} \cap  Ker d^{\prime}
$$
and
$$
\frac
{d^{\prime}F^{p-1,q} \cap  Ker d^{\prime\prime} \cap  F^{p,q}}
{d^{\prime}F^{p-1,q} \cap  d^{\prime\prime}F^{p,q-1}}
\backsimeq
\bold {H_{d^{\prime\prime}}} \cap  d^{\prime}F^{p-1,q}.
$$
\par
$\bold {key equality}$.
$$
d^{\prime}\square _{d^{\prime\prime}}u
=\frac {p+q-n+2}{p+q-n+1}d^{\prime\prime}d^{\prime\prime\ast}(d^{\prime}u)
+\frac {p+q-n+3}{p+q-n+2}d^{\prime\prime\ast}d^{\prime\prime}(d^{\prime}u)
$$
$$
for u \in  F^{p,q}.
$$
We show this.
$$
\aligned
d&^{\prime}\{d^{\prime\prime}(\delta ^{\prime\prime}u-
\frac 1{p+q-n+1}\Lambda L \delta ^{\prime\prime}u)\}\\
=&d^{\prime\prime}\{-d^{\prime}(\delta ^{\prime\prime}u
-\frac 1{p+q-n+1}\Lambda L \delta ^{\prime\prime}u)\}\\
=&d^{\prime\prime}\{\delta ^{\prime\prime}d^{\prime}u
+\frac 1{p+q-n+1}d^{\prime}\Lambda L \delta ^{\prime\prime}u\}
\endaligned
$$
While
$$
d^{\prime}\Lambda - \Lambda d^{\prime}=-\sqrt {-1}\delta ^{\prime\prime}.
$$
Therefore
$$
\aligned
d^{\prime}\Lambda L\delta ^{\prime\prime}u=&
\Lambda d^{\prime}L\delta ^{\prime\prime}u-
\sqrt {-1}\delta ^{\prime\prime}L\delta ^{\prime\prime}u ( by 
\delta ^{\prime\prime}L-L\delta ^{\prime\prime}=-\sqrt {-1}d^{\prime})\\
=&-\sqrt {-1}\delta ^{\prime\prime}L\delta ^{\prime\prime}u\\
=&\delta ^{\prime\prime}d^{\prime}.
\endaligned
$$
So we have
$$
\aligned
d^{\prime}d^{\prime\prime}d^{\prime\prime\ast}u
=&\frac {p+q-n+2}{p+q-n+1}d^{\prime\prime}
\delta ^{\prime\prime\ast}d^{\prime}u\\
=&\frac {p+q-n+2}{p+q-n+1}d^{\prime\prime}
(\delta ^{\prime\prime}d^{\prime}u-\frac 1{p+q-n+2}\Lambda L 
\delta ^{\prime\prime}d^{\prime}u)\\
=&\frac {p+q-n+2}{p+q-n+1}d^{\prime\prime}d^{\prime\prime\ast}d^{\prime}u.
\endaligned
$$
By the same way, we have
$$
d^{\prime}d^{\prime\prime\ast}d^{\prime\prime}u
=\frac {p+q-n+3}{p+q-n+2}d^{\prime\prime\ast}d^{\prime\prime}d^{\prime}u.
$$
 The correspondence is as follows.
$$
\frac {Ker d' \cap   Ker d^{\prime\prime} \cap  F^{p,q}}{Ker d' \cap
d^{\prime\prime}F^{p,q-1}}  \to  \bold H_{d^{\prime\prime}} \cap  Ker d'
$$
for $u  \to  \bold H_{d^{\prime\prime}}u$.\par
And
$$
\frac {d'F^{p-1,q}  \cap  Ker d^{\prime\prime} \cap  F^{p,q}}
{d'F^{p-1,q} \cap  d^{\prime\prime}F^{p,q-1}}  \to 
\bold H_{d^{\prime\prime}} \cap  d'F^{p-1,q}
$$
for $u  \to  \bold H_{d^{\prime\prime}}u$.\par

 $\bold {Lemma 5.2}$. If $p+q\ge n+1$, then for $u$ in $F^{p,q}$,
$$
d'\bold H_{d^{\prime\prime}}u=\bold H_{d^{\prime\prime}}d'u .
$$
$\bold {Proof}$. For $u$ in $F^{p,q}$,
$$
d'u=\bold H_{d^{\prime\prime}}d'u+
\square _{d^{\prime\prime}}N_{d^{\prime\prime}}d'u.
$$
And
$$
u=\bold H_{d^{\prime\prime}}u+
\square _{d^{\prime\prime}}N_{d^{\prime\prime}}u.
$$
While by the key equality,
$$
d^{\prime}\square _{d^{\prime\prime}}u
=\frac {p-1+q-n+2}{p-1+q-n+1}d^{\prime\prime}d^{\prime\prime\ast}(d^{\prime}u)
+\frac {p-1+q-n+3}{p-1+q-n+2}d^{\prime\prime\ast}d^{\prime\prime}(d^{\prime}u)
$$
We take $v=\bold H_{d^{\prime\prime}}u$ and put $v$ in the place of $u$, then
$$
\aligned
d'\square _{d^{\prime\prime}}\bold H_{d^{\prime\prime}}u
=&\frac {p-1+q-n+2}{p-1+q-n+1}d^{\prime\prime}d^{\prime\prime\ast}(d^{\prime}u)
+\frac {p-1+q-n+3}{p-1+q-n+2}d^{\prime\prime\ast}d^{\prime\prime}(d^{\prime}u)
\\
0=&\frac {p-1+q-n+2}{p-1+q-n+1}d^{\prime\prime}d^{\prime\prime\ast}(d^{\prime}u)
+\frac {p-1+q-n+3}{p-1+q-n+2}d^{\prime\prime\ast}d^{\prime\prime}(d^{\prime}u)
\endaligned
$$
Therefore 
$d'\bold H_{d^{\prime\prime}}u$ is a harmonic form.
This means that
$$
\bold H_{d^{\prime\prime}}d'u=d'\bold H_{d^{\prime\prime}}u
$$
by taking the harmonic part.\par
$\bold {Q.E.D.}$\par

Even for the case $p+q=n$, we can introduce a harmonic operator by;
$$
for u in F^{p,q},
$$
$$
\bold H_{d^{\prime\prime}}u = u -
d^{\prime\prime\ast}N_{d^{\prime\prime}}d^{\prime\prime}u .
$$

\enddocument